\documentclass[reqno]{amsart}
\usepackage{hyperref}

\def\N{{\mathbb{N}}}
\def\Q{{\mathbb{Q}}}
\def\R{{\mathbb{R}}}

\begin{document}
\title[continuity]{On the almost everywhere continuity}

\author[Blot]
{Jo${\rm \ddot e}$l Blot}

\address{Jo\"{e}l Blot: Laboratoire SAMM EA 4543,\newline
Universit\'{e} Paris 1 Panth\'{e}on-Sorbonne, centre P.M.F.,\newline
90 rue de Tolbiac, 75634 Paris cedex 13,
France.}
\email{blot@univ-paris1.fr}
\date{November 4, 2014}

\numberwithin{equation}{section}
\newtheorem{theorem}{Theorem}[section]
\newtheorem{lemma}[theorem]{Lemma}
\newtheorem{example}[theorem]{Example}
\newtheorem{remark}[theorem]{Remark}
\newtheorem{definition}[theorem]{Definition}
\newtheorem{proposition}[theorem]{Proposition}
\newtheorem{corollary}[theorem]{Corollary}
\begin{abstract}
The aim of this paper is to provide characterizations of the Lebesgue-almost everywhere continuity of a function $f : [a,b] \rightarrow \R$. These characterizations permit to obtain necessary and sufficient conditions for the Riemann integrability of $f$.
\end{abstract}
\maketitle
M.S.C. 2010: 26A15, 26A42.
\section{Introduction}
The main aim of this paper is to establish the following theorem.
\begin{theorem}\label{th11}
Let $a >b$ be two real numbers, and $f : [a,b] \rightarrow \R$ be a function. We assume that $f$ admits a finite right-hand limit at each point of $[a,b)$ except on a Lebesgue-negligible set (respectively on a at most countable set). Then $f$ is continuous at each point of $[a,b]$ except on a Lebesgue-negligible set (respectively on a at most countable set).
\end{theorem}
\vskip2mm
The origine of this work is a paper of Daniel {\sc Saada} \cite{Sa} which states that a real function defined on a real segment which right-hand continuous possesses at most an at most countable subset of discontinuity points. Saada attributes the proof of this result to Alain {\sc R\'emondi\`ere}. Studying this result and its proof, we see that it contains a central argument that we have described in our Lemma \ref{lem41} and Lemma \ref{lem42}, and we use this argument to obtain other results. And so the present work is a continuation of the work of R\'emondi\`ere and Saada. 
\vskip2mm
In Section 2 we precise our notation and we give comments on them. In Section 3, we establish lemmas which are useful for the proof of Theorem \ref{th11}. In Section 4, we provide results on the left-hand continuity and on the right-hand continuity. In Section 5 we give the proof of Theorem \ref{th11}. In Section 6 we establish corollaries of Theorem \ref{th11}.
\section{Notation}
We use the left-hand oscillation of $f$ at $x \in (a,b]$ defined by
$$\omega_L(x) := \lim_{h \rightarrow 0+}( \sup_{y \in [x-h,x]} f(y)) - \lim_{h \rightarrow 0+}( \inf_{y \in [x-h,x]} f(y))$$
and also the righ-hand oscillation of $f$ at $x \in [a,b)$ defined by 
$$\omega_R(x) := \lim_{h \rightarrow 0+}( \sup_{y \in [x,x+h]} f(y)) - \lim_{h \rightarrow 0+}( \inf_{y \in [x,x+h]} f(y)).$$
Note that we have $\sup_{y \in [x-h,x]} f(y) \geq f(x)$ and consequently 
$$\lim_{h \rightarrow 0+}( \sup_{y \in [x-h,x]} f(y)) \geq f(x) > - \infty.$$
 Also note that we have 
$ \inf_{y \in [x-h,x]} f(y) \leq f(x)$ and consequently 
$$\lim_{h \rightarrow 0+}( \inf_{y \in [x-h,x]} f(y)) \leq f(x) < + \infty.$$
 And so, $\omega_L(x)$ is a sum of two elements of $(- \infty, + \infty]$ and therefore it is well-defined in  $(- \infty, + \infty]$; more precisely its belongs to $[0, + \infty]$. For similar reasons, $\omega_R(x)$ is well-defined in  $[0, + \infty]$.
\vskip2mm
We use the following notation when $g : [a,b] \rightarrow (-\infty, + \infty]$ and $r \in \R$:\\
 $\{g=0 \} := \{x \in [a,b] : g(x) = 0 \}$, $\{g >r \} := \{x \in [a,b] : g(x) > r \}$,\\
   $\{g \leq r \} := \{x \in [a,b] : g(x) \leq  r \}$.
\vskip2mm
A subset $N \subset [a,b]$ is called {\it Lebesgue-negligible} when there exists $B$, a borelian subset of $[a,b]$, such that $N \subset B$ and $\mu(B) = 0$ where $\mu$ denotes the Lebesgue measure of $\R$. Such a vocabulary is used for instance in \cite{Ch1}.
\vskip2mm
\begin{remark}\label{rem21}
The following equivalence hold.
\begin{enumerate}
\item[(A)] When $x \in (a,b]$, $\omega_L(x) = 0$ if and only if $f$ is left-hand continuous at $x$.
\item[(B)] When $x \in [a,b)$, $\omega_R(x) = 0$ if and only if $f$ is right-hand continuous at $x$.
\item[(C)] When $x \in (a,b)$, ($\omega_L(x) = 0$ and $\omega_R(x) = 0$)  if and only if $f$ is continuous at $x$.
\end{enumerate}
\end{remark}
These equivalences are easy to prove. One important fact is that $x$ belongs to the neighborhoods $[x-h,x]$ and $[x, x+h]$.
\begin{remark}\label{rem22} When we will speak of the left-hand limit (respectively of the right-hand limit) of the function $f$ at $x$, we speak of the limit of $f(y)$ when
$y \rightarrow x, y>x$ (respectively $y \rightarrow x, y< x$); the point $x$ is not included into his "`neighborhoods"'. The situation is different in the definition of the oscillations
$\omega_L$ and $\omega_R$. We denote $f(x-) := \lim\limits_{y \rightarrow x, y<x} f(y)$ and $f(x+) := \lim\limits_{y \rightarrow x, y>x} f(y)$.
\end{remark}

\section{preliminaries}
We establish lemmas which are useful to the proof of Theorem \ref{th11}. 
\begin{lemma}\label{lem31}
Let $f : [a,b] \rightarrow \R$ be a function, and $z \in [a,b)$. We assume that $f$ admits a finite right-hand limit at $z$. Then we have :\\
$$\forall \epsilon > 0, \exists \lambda(z,\epsilon) > 0, \forall x \in (z, z +  \lambda(z,\epsilon)], \omega_L(x) \leq \epsilon.$$
\end{lemma}
\vskip4mm
\begin{proof}
We arbitrarily fix $\epsilon > 0$. Using the assumption, there exists $d_z \in \R$ such that
\begin{equation}\label{eq31}
\exists \eta(z,\epsilon) > 0, \forall x \in [a,b], z < x \leq z + \eta(z,\epsilon) \Longrightarrow \vert f(x) - d_z \vert \leq \epsilon.
\end{equation}
When $x \in (z, z + \eta(z, \frac{\epsilon}{4})]$ and when $y \in (z,x]$, we have $y \in (z, z + \eta(z, \frac{\epsilon}{4})]$, and using (\ref{eq31}) we obtain $
\vert f(x) - d_z \vert \leq \frac{\epsilon}{4}$ and $\vert f(x) - d_z \vert \leq \frac{\epsilon}{4}$ that implies
$$\vert f(x) - f(y) \vert \leq \vert f(x) - d_z \vert + \vert f(y) - d_z \vert \leq 2 \frac{\epsilon}{4} = \frac{\epsilon}{2} \Longrightarrow$$
$$f(x) - \frac{\epsilon}{2} \leq f(y) \leq f(x) + \frac{\epsilon}{2}.$$
Then, for all $h \in (0, x-z]$, we obtain
$$f(x) - \frac{\epsilon}{2} \leq \inf_{y \in [x-h,x]} f(y) \leq \sup_{y \in [x-h,x]} f(y)) \leq f(x) + \frac{\epsilon}{2} \Longrightarrow$$
$$f(x) - \frac{\epsilon}{2} \leq \lim_{h \rightarrow 0+}(\inf_{y \in [x-h,x]} f(y)) \leq \lim_{h \rightarrow 0+}(\sup_{y \in [x-h,x]} f(y)) \leq f(x) + \frac{\epsilon}{2}$$
$$\Longrightarrow 0 \leq \omega_L(x) \leq f(x) + \frac{\epsilon}{2} -(f(x) - \frac{\epsilon}{2}=  \epsilon.$$
And so it suffices to take $\lambda(z, \epsilon) := \eta(z, \frac{\epsilon}{4})$.
\end{proof}
Using similar arguments we can prove the following result.
\vskip2mm
\begin{lemma}\label{lem32}
Let $f : [a,b] \rightarrow \R$ be a function, and $z \in (a,b]$. We assume that $f$ admits a finite left-hand limit at $z$. Then we have :\\
$$\forall \epsilon > 0, \exists \nu(z,\epsilon) > 0, \forall x \in [z - \nu(z, \epsilon), z ), \omega_R(x) \leq \epsilon.$$
\end{lemma}
\vskip2mm
\begin{lemma}\label{lem33}
Let $I$ be a nonempty set, and $(S_i)_{i \in I}$ be a family of subintervals of $[a,b]$ such that $S_i \cap S_j = \emptyset$ when $i \neq j$, and such $\mu(S_i) >0$ for all $i \in I$, where $\mu$ denotes the Lebesgue measure of $\R$. Then $I$ is at most countable.
\end{lemma}
\vskip2mm
\begin{proof}
Since a positive measure is additive, for all finie subset $J \subset I$, we have $\mu (\cup_{j \in J} S_j) = \sum_{j \in J} \mu(S_j)$. Since a positive measure is monotonic, $
\cup_{j \in J} S_j \subset [a,b]$ implies $\mu(\cup_{j \in J} S_j) \leq \mu([a,b])= b-a$, and so we have $\sum_{j \in J} \mu(S_j) \leq b-a < + \infty$ for all finite subset $J$ of $I$. Therefore the family of non negative real numbers $(\mu(S_i))_{i \in I}$ is summable in $[0, + \infty)$, and consequently the set $\{ i \in I : \mu(S_i) \neq 0 \}$ is at most countable (Corrolary 9-9, p. 220 in \cite{Ch2}). Since $\mu(S_i) >0$ for all $i \in I$, we obtain that $I$ is at mot countable.
\end{proof}
\begin{remark}\label{rem34}
We can also prove Lemma \ref{lem33} by building a function $\varphi : I \rightarrow \Q$ in the following way: since $\Q$ is dense into $\R$, for each $i \in I$, there exists $\varphi (i) \in \Q \cap S_i$. Since $I_i \cap I_j = \emptyset $ when $i \neq j$, we have $\varphi (i) \neq \varphi (j)$ when $i \neq j$. And so $\varphi$ is injective. Since $\Q$ is countable, $\varphi(I) \subset \Q$ is at most countable, and using an abridgement of $\varphi$, we build a bijection between $\varphi(I)$ and $I$.
\end{remark}
\section{Limits on one side, continuities on the other side}
The following results establish that the existence of left-hand (respectively right-hand) limits implies the right-hand (respectively left-hand) continuity.
\vskip2mm
\begin{lemma}\label{lem41}
let $f : [a,b] \rightarrow \R$ be a function, and $N$ be a Lebesgue-negligible (respectively at most countable) subset of $[a,b]$. We assume that $f$ admits a finite right-hand limit at each $x \in [a,b) \setminus N$.\\
Then the set of the points of $[a,b]$ where f is not left-hand continuous is Lebesgue-negligible (respectively at most countable).
\end{lemma}
\vskip4mm
\begin{proof}
We arbitrarily fix $\epsilon > 0$. Using Lemma \ref{lem31}, denoting $\lambda_z := \lambda(z, \epsilon)$, we obtain the following assertion.
\begin{equation}\label{eq41}
\forall z \in \{\omega_L > \epsilon \} \cap ([a,b] \setminus N), \exists \lambda_z > 0, (z, z + \lambda_z] \subset \{ \omega_L \leq \epsilon \}.
\end{equation}
Let $z_1, z_2 \in \{\omega_L > \epsilon \} \cap ([a,b] \setminus N)$, $z_1 \neq z_2$. We can assume that $z_1 < z_2$. After (\ref{eq41}), we cannot have $z_2$ into 
$(z_1, z_1 + \lambda_{z_1}]$, therefore we have $ z_2 > z_1 + \lambda_{z_1}$, and we have proven:
$$\forall z_1, z_2 \in \{\omega_L > \epsilon \} \cap ([a,b] \setminus N), z_1 \neq z_2 \Longrightarrow (z_1, z_1 + \lambda_{z_1}] \cap (z_2, z_2 + \lambda_{z_2}] = \emptyset.$$
We have also $\mu((z, z + \lambda_z]) =  \lambda_z > 0$. Then using Lemma \ref{lem33}, we can assert that 
\begin{equation}\label{eq42}
\forall \epsilon > 0,  \{\omega_L > \epsilon \} \cap ([a,b] \setminus N) \; \; {\rm is} \; {\rm at} \; {\rm most} \; {\rm countable}.
\end{equation}
Note that $\{ \omega_L > 0 \} = \omega_L^{-1}((0, +\infty]) = \omega_L^{-1}(\bigcup_{n \in \N_*} (\frac{1}{n}, + \infty])$\\ = $\bigcup_{n \in \N_*}\omega_L^{-1}((\frac{1}{n}, + \infty]) = 
 \bigcup_{n \in \N_*} \{ \omega_L > \frac{1}{n} \} \Longrightarrow$
$$\{ \omega_L > 0 \} \cap  ([a,b] \setminus N) = \bigcup_{n \in \N_*} (\{ \omega_L > \frac{1}{n} \} \cap ([a,b] \setminus N)).$$
Using (\ref{eq42}), since a countable union of at most countable subsets is at most countable, we ontain the following assertion.
\begin{equation}\label{eq43}
  \{\omega_L > 0 \} \cap ([a,b] \setminus N) \; \; {\rm is} \; {\rm at} \; {\rm most} \; {\rm countable}.
\end{equation}
Note that $\{\omega_L > 0 \} = (\{\omega_L > 0 \} \cap ([a,b] \setminus N)) \cup (\{\omega_L > 0 \} \cap N)$.\\
 Since $(\{\omega_L > 0 \} \cap N) \subset N$ and since $N$ is Lebesgue-negligible (respectively at most countable), $(\{\omega_L > 0 \} \cap N)$ is Lebesgue-negligible (respectively at most countable). Recall that an at most countable subset of $\R$ is Lebesgue-negligible. And so when $N$ is Lebesgue-negligible,  $\{\omega_L > 0 \}$ is Lebesgue-negligible as a union of two Lebesgue-negligible susbets, and when $N$ is at most countable, $\{\omega_L > 0 \}$ is at most countable as a union of two at most countable subsets. Using (A) of Remark \ref{rem21}, the lemma is proven.
\end{proof}
Proceegings as in the proof of Lemma \ref{lem41}, we obtain the following result.
\vskip2mm
\begin{lemma}\label{lem42}
let $f : [a,b] \rightarrow \R$ be a function, and $M$ be a Lebesgue-negligible (respectively at most countable) subset of $[a,b]$. We assume that $f$ admits a finite left-hand limit at each $x \in (a,b] \setminus M$.\\
Then the set of the points of $[a,b]$ where f is not right-hand continuous is Lebesgue-negligible (respectively at most countable).
\end{lemma}
\section{Proof of Theorem \ref{th11}}
Using Lemma \ref{lem41} and (A) of Remark \ref{rem21}, $\{ \omega_L > 0 \}$ is Lebesgue-negligible (respectively at most countable) since $\{ \omega_L > 0 \}$ is exactly the set of the points of $(a,b]$ where $f$ is not left-hand continuous.
\vskip1mm
Now, setting $M = \{ \omega_L > 0 \}$, for all $x \in [a,b] \setminus M$, $f(x-) = f(x) \in \R$, and the assumption of Lemma \ref{lem42} is fulfilled. Consequently we obtain that 
 $\{ \omega_R > 0 \}$ is Lebesgue-negligible (respectively at most countable) after (B) of Remark \ref{rem21}.
\vskip1mm
Note that $\{ \omega_L = 0 \} \cap \{ \omega_R = 0 \}$ is exactly the set of the points of $(a,b)$ where $f$ is continuous. We have 
$$[a,b] \setminus (\{ \omega_L = 0 \} \cap  \{ \omega_R = 0 \}) = [a,b]  \cap (\{ \omega_L > 0 \} \cup  \{ \omega_R > 0 \}) = \{ \omega_L > 0 \} \cup  \{ \omega_R > 0 \}.$$
This set is Lebesgue-negligible (respectively at most countable) as a union of two Lebesgue-negligible (respectively at most countable) sets. Note that $\{ a, b \}$ is Lebesgue-negligible (respectively at most countable) and so set of the discontinuity points of $f$ is Lebesgue-negligible (respectively at most countable).
\section{Consequences}
A first consequence of Theorem \ref{th11} is the following result.
\vskip2mm
\begin{theorem}\label{th61}
Let $a >b$ be two real numbers, and $f : [a,b] \rightarrow \R$ be a function. Then the following assertions are equivalent.
\begin{enumerate}
\item[($\alpha$)] The set of the discontinuity points of $f$ is Lebesgue-negligible (respectively at most countable).
\item[($\beta$)] The set of the left-hand discontinuity points of $f$ is Lebesgue-negligible (respectively at most countable).
\item[($\gamma$)] The set of the right-hand discontinuity points of $f$ is Lebesgue-negligible (respectively at most countable).
\item[($\delta$)] The set of the points where $f$ does not admit a finite left-hand limit is Lebesgue-negligible (respectively at most countable).
\item[($\epsilon$)] The set of the points where $f$ does not admit a finite right-hand limit is Lebesgue-negligible (respectively at most countable).
\end{enumerate}
\end{theorem}
\vskip2mm
\begin{proof}
The implications $(\alpha) \Longrightarrow (\beta) \Longrightarrow (\delta)$ are easy, and $ (\delta) \Longrightarrow (\alpha)$ is Theorem \ref{th11}. The implications $
(\alpha) \Longrightarrow (\gamma)  \Longrightarrow (\epsilon)$ are easy. we can do a proof which is similar to this one of Theorem \ref{th11} to prove $ (\epsilon) \Longrightarrow (\alpha)$.
\end{proof}
About the Riemann-integrability we recall a famous theorem of Lebesgue, \cite{Le} p. 29, \cite{SG} p. 20.
\vskip2mm
\begin{theorem}\label{th62}
Let $a >b$ be two real numbers, and let $f : [a,b] \rightarrow \R$ be a bounded function. Then the following assertions are equivalent.
\begin{enumerate}
\item[(i)] $f$ is Riemann integrable on $[a,b]$.
\item[(ii)] The set of the discontinuity points of $f$ is Lebesgue-negligible.
\end{enumerate}
\end{theorem}
\vskip3mm
As a  consequence of Theorem  \ref{th61} and of the previous classical theorem of Lebesgue, we obtain the following result on the Riemann-integrability.
\vskip2mm
\begin{theorem}\label{th63}
Let $a >b$ be two real numbers, and let $f : [a,b] \rightarrow \R$ be a bounded function. Then the following assertions are equivalent.
\begin{enumerate}
\item[({\sf a})] $f$ is Riemann integrable on $[a,b]$.
\item[({\sf b})] The set of the points where $f$ does not admit a finite left-hand limit is Lebesgue-negligible.
\item[({\sf c})] The set of the points where $f$ does not admit a finite right-hand limit is Lebesgue-negligible.
\end{enumerate}
\end{theorem}
An easy consequence of this result is the following one.
\begin{corollary}\label{cor64}
Let $a >b$ be two real numbers, and let $f : [a,b] \rightarrow \R$ be a function. If $f$ is right-hand continuous on $[a,b]$ or left-hand continuous on $[a,b]$, then the set of the discontinuity points of $f$ is at most countable, and consequently when in addition $f$ is assumed to be bounded, $f$ is Riemann integrable on $[a,b]$.
\end{corollary}
{\bf Acknowledgements.} I thanks my colleagues B. Nazaret, M. Bachir and J.-B. Baillon for interesting discussions on these topics.

\end{document}